\documentclass[12pt]{article}
\usepackage{amsfonts,amssymb,amsmath,amsthm}
\usepackage{graphicx}
\usepackage{epstopdf}
\ifpdf
  \DeclareGraphicsExtensions{.eps,.pdf,.png,.jpg}
\else
  \DeclareGraphicsExtensions{.eps}
\fi
\usepackage{pbox}
\usepackage{color,url}
\usepackage{dsfont}
\usepackage{setspace}
\usepackage{graphicx}
\usepackage{epsfig}
\usepackage{float}
\usepackage{stmaryrd}
\usepackage{pifont}
\usepackage{comment}
\usepackage{enumerate}
\usepackage[section]{placeins}

\newtheorem{thm}{\bf{Theorem}}[section]
\newtheorem{lem}[thm]{\bf{Lemma}}

\newtheorem{df}[thm]{\bf{Definition}}
\newtheorem{cor}[thm]{\bf{Corollary}}

\newtheorem{ex}[thm]{\bf{Example}}

\newcommand{\dom}{\operatorname{dom}}

\newcommand{\1}{\operatorname{\mathds{1}}}

\newcommand{\Id}{\operatorname{Id}}

\newcommand{\diag}{\operatorname{diag}}
\newcommand{\Diag}{\operatorname{Diag}}

\newcommand{\R}{\operatorname{\mathbb{R}}}

\newcommand{\bbm}{\begin{bmatrix}}
\newcommand{\ebm}{\end{bmatrix}}

\newcommand{\one}{\mathbf{1}}

\newcommand{\cshd}{\diag \left [\nabla^2_c f(x^0;S) \right ]}
\newcommand{\diagth}{\diag [\nabla^2f(x^0) ]}
\newcommand{\Diagth}{\Diag [\nabla^2f(x^0) ]}
\newcommand{\upperh}{U[\nabla^2 f(x^0)]}
\newcommand{\offdiag}{N[\nabla^2 f(x^0)]}

\DeclareMathOperator{\RMPB}{RMPB}
\DeclareMathOperator{\CMPB}{CMPB}
\DeclareMathOperator{\CB}{CB}
\DeclareMathOperator{\RB}{RB}

\DeclareMathOperator{\RER}{RE}

\DeclareMathOperator{\POI}{POI}

\begin{document}
\title{Approximating  the diagonal of a Hessian: which sample set of points should be used}
\author{Gabriel Jarry--Bolduc\thanks{Department of Mathematics, University of British Columbia, Okanagan Campus, Kelowna, B.C. V1V 1V7, Canada. Research partially supported by Natural Sciences and Engineering Research Council (NSERC) of Canada Discovery Grant 2018-03865. gabjarry@alumni.ubc.ca }}
\maketitle\author

\begin{abstract}
An explicit formula to approximate the diagonal entries of the Hessian is introduced. When the derivative-free technique called \emph{generalized centered simplex gradient} is used to approximate the gradient, then the formula can be computed for only one additional function evaluation. An error bound is introduced and provides information on  the form of the sample set of points that should be used to approximate the diagonal of a Hessian. If the sample set of points is built in a specific manner, it is shown that the technique is $\mathcal{O}(\Delta_S^2)$ accurate approximation of the diagonal entries of the Hessian where $\Delta_S$ is the radius of the sample set. 
\end{abstract}

\noindent{\bf Keywords:} (generalized) centered simplex gradient; centered simplex Hessian diagonal; derivative-free optimization.


\section{Introduction}
Derivative-free optimization (DFO) constructs optimization algorithms that do not employ first-order information within the algorithm. Recently,  substantial progress has been made regarding their applications and numerical implementations (see \cite{Audet2014,Audet2017,Conn2009,Custodio2017,HareNutiniTesfamariam2013,larson_menickelly_wild_2019}).

One of the main categories of DFO algorithms is model-based DFO methods.  
A simple method to approximate the objective function is to build a linear interpolation model  from $n+1$ well-poised sample points in $\R^n$.  The gradient of this linear model is  called the {\em simplex gradient} and provides an approximation of the true gradient \cite{Bortz1998,Kelley1999}. 

An error bound comparing the simplex gradient and the true gradient was introduced in  the late 1990s and it is order $\mathcal{O}(\Delta_S)$, where $\Delta_S$ is the radius of the  sample set of evaluated points \cite{Kelley1999}. This error bound shows that the optimizer can control the accuracy of the approximation  technique by varying the radius $\Delta_S$ of the sample set of points
\cite[Ch.\ 10 \& 11]{Audet2017}. A generalization of the simplex gradient called \emph{generalized simplex gradient} has the advantage of not being limited to the setting where exactly $n+1$ interpolation points are used in $\R^n$. In \cite{Conn22008}, the authors consider the case where less than $n+1$ sample points (underdetermined models) and more  than $n+1$ points (overdetermined models)  are used to approximate the gradient.  Most importantly, they establish an error bound for the overdetermined case  and show that  it retains order  $\mathcal{O}(\Delta_S).$  This topic is also  investigated in \cite{Regis2015} and calculus rules for the generalized simplex gradient  are introduced  in \cite{Regis2015, hare2020}.

Many other methods of approximating gradients exist \cite{Billups2013,Oeuvray2007,Powell2004,Regis2005,Schonlau1998,Wild2011}. One of these methods is the \emph{centered simplex gradient} \cite{Kelley1999}. Using the properties of the Moore-Penrose pseudoinverse (see Definition \ref{def:mpinverse}), the approach is generalized so that it does not require  exactly $2n$  sample points in $\R^n.$ The approach is called {\em generalized centered simplex gradient} (GCSG) and it is created by retaining the $k$ original points in the sample set and adding their reflection through the point of interest (see Definition \ref{def:cshd}).  An error bound which applies to the underdetermined, determined and overdetermined cases is introduced in \cite{hare2020error}. The error bound shows  that  the  GCSG is $\mathcal{O}(\Delta_S^2)$  accurate. 

In 2020,  Tappenden and Coope  showed how to compute the GCSG in $\mathcal{O}(n)$ flops when using four different sample sets of points \cite{coope2021gradient}. For these specific choices of sample sets, a formula to approximate the diagonal entries of the Hessian is provided. We can observe that if the gradient is approximated via the  GCSG, then only one additional function evaluation is sufficient to obtain an approximation of the diagonal entries of the Hessian. Numerical examples are provided to compare the accuracy of the gradient and diagonal entries of the Hessian depending on  the sample set of points utilized. These numerical examples  agree with the error bound defined for the gradient, but no error bound for the diagonal entries of the Hessian is provided. For this reason, the poor accuracy obtained using certain sample set of points is not fully  explained. In this paper, an error bound for the diagonal entries of the Hessian is introduced. Analyzing this error bound, we will obtain valuable information on the form  of the sample set that should be utilized. In particular, this will explain the poor accuracy obtained with some of the sample sets of points used in the numerical examples in \cite{coope2021gradient}. Moreover,  an explicit formula to compute the diagonal entries of a Hessian for any  nonempty sample set of points  will be introduced. The technique is called \emph{centered simplex Hessian diagonal} (CSHD).

The structure of this paper is the following. In Section \ref{sec:prel}, we introduce notation and basic definitions. In Section \ref{sec:cshd}, the formula for the CSHD is introduced and an error bound is proven. When the sample set of point has a specific form, it is shown that the CSHD is  $\mathcal{O}(\Delta_S^2)$ accurate. In Section \ref{sec:num}, numerical examples are provided for different sample sets of points and an error analysis is provided.  Section \ref{sec:conc} summarizes the work accomplished and suggests some topics to explore in future research.

\section{Preliminaries} \label{sec:prel}
Unless otherwise stated, we use the standard notation found in \cite{rockwets}. The domain of a function $f$ is denoted by $\dom f$. The transpose of a matrix $A$ is denoted by $A^\top$. We work in finite-dimensional space $\R^n$ with inner product $x^\top y=\sum_{i=1}^nx_iy_i$. The norm of a vector is denoted $\Vert x \Vert$ and is taken to be the $\ell_2$ norm. 
Given a matrix $A \in \R^{n \times k},$ the  $\ell_2$ induced matrix norm is used. That is 
\begin{align*}
    \left \Vert A \right \Vert_2&=\left \Vert A \right \Vert= \max \{\Vert Ax \Vert: \Vert x \Vert=1 \}.
\end{align*}
We denote by $B(x^0,\Delta)$ the close ball centered about $x^0$ with radius $\Delta$.
  Let $M \in \R^{n \times n}$. We will use the notation $\diag [M]$  to denote the vector in $\R^n$ containing the diagonal entries of $M$. The notation $\Diag[M] \in \R^{n \times n}$  represents the matrix containing the diagonal entries of $M.$ We will denote by $U[M] \in \R^{n \times n}$ the strictly upper triangular entries of $A$ and by $N[M] \in \R^{n \times n},$ the matrix containing the off-diagonal entries of a matrix $M$.  The identity matrix in $\R^{n \times n}$ is denoted by $\Id_n$ and the vector of all ones in $\R^n$ is denoted by $\one_n.$
 
 Recall that a generalization of the matrix inverse is the Moore--Penrose pseudoinverse.
\begin{df}[Moore--Penrose pseudoinverse] \cite{Horn1990} \label{def:mpinverse}
Let $A\in\R^{n\times k}$. The \emph{Moore--Penrose pseudoinverse} of $A$, denoted by $A^\dagger$, is the unique matrix in $\R^{k \times n}$ that satisfies the following four equations:
\begin{equation*}
(1) \, AA^\dagger A=A, \quad (2) \, A^\dagger AA^\dagger=A^\dagger,  \quad (3) \, (AA^\dagger)^\top=AA^\dagger, \quad (4) \, (A^\dagger A)^\top=A^\dagger A.
\end{equation*}
\end{df}
The Moore--Penrose inverse $A^\dagger$ satisfies  the following three properties.
\begin{enumerate}[(1)]
\item If $A$ has full column rank $k$, then $A^\dagger$ is a left-inverse of $A$, that is, $A^\dagger A=\Id_k$.
\item If $A$ has full row rank $n$, then $A^\dagger$ is a right-inverse of $A$, that is, $AA^\dagger=\Id_n$.
\item If $A$ is a square matrix with full rank, then $A^\dagger=A^{-1},$ the inverse of $A.$
\end{enumerate}
 Last, recall the definition of the \emph{Hadamard product}.

\begin{df}\cite{Horn1990}
Let $A=\bbm a_{i,j} \ebm \in \R^{n \times k}$ and $B=\bbm b_{i,j}\ebm \in \R^{n \times k}.$ The Hadamard product of $A$ and $B$, denoted $A \odot B$ is the component wise product. That is $A \odot B=\bbm a_{i,j} b_{i,j} \ebm \in \R^{n \times k}.$
\end{df}
\section{The CSHD and its error bound} \label{sec:cshd}
We now present a process to obtain a formula  to approximate the diagonal entries of the Hessian. This process is assuming that the gradient is approximated via the GCSG \cite{hare2020error}.
 Let $f:\dom f \subseteq \R^n \to \R$ be $\mathcal{C}^2$, $x^0 \in \dom f$ be the point of interest, and  $S=\bbm s^1&s^2& \cdots &s^k\ebm \in \R^{n \times k}$ be a set of nonzero distinct directions in $\R^n$ written in matrix form. The matrix $S$ contains all  the directions to add and subtract to the point of interest $x^0$ to form the sample points where the function $f$ is evaluated. Suppose a model $m(x)$ is constructed such that
 $$f(x^0+s^i)=m(x^0+s^i), \quad i \in \{1, 2, \dots, k\},$$ where $m$ is a diagonal quadratic model. That is  
 \begin{equation}\label{eq:model}
 m(x)= f(x^0)+ g^\top (x-x^0)+\frac{1}{2}(x-x^0)^\top D (x-x^0)
 \end{equation}
 where $D \approx \Diag \left [\nabla^2 f(x^0) \right ] \in \R^{n \times n}$ and $g \approx \nabla f(x^0) \in \R^n.$ Define $$ W= \bbm s^1 \odot s^1 &\cdots &s^k \odot s^k \ebm \in \R^{n \times k},  \, \text{and} \quad 
  \delta_f (x^0;S)=\bbm f(x^0+s^1)-f(x^0)\\ \vdots \\ f(x^0+s^k)-f(x^0)\ebm \in \R^k.$$
 Letting  $x=x^0+s^i,i \in  \{1, 2, \dots, k\}$, in \eqref{eq:model}, we get the system
 \begin{equation} \label{eq:system}
     \delta_f(x^0;S)=S^\top g +\frac{1}{2} W^\top d
 \end{equation}
 where $d=\diag \left [D(x^0) \right ] \in \R^n.$
  When the GCSG is employed, the sample points $x^0-s^i$ for all $i \in \{1, \dots, k\} $ are also created. Letting $x=x^0-s^i$ for all $i$ in Equation \eqref{eq:model}, we obtain the system
 \begin{equation} \label{eq:system2}
     \delta_f(x^0;-S)=-S^\top g+\frac{1}{2}  W^\top d.
 \end{equation}
 System \eqref{eq:system} and System \eqref{eq:system2} can be combined to form one  block matrix equation:
 \begin{equation} \label{eq:blocksystem}
 \bbm \delta_f(x^0;S)\\ \delta_f (x^0;-S) \ebm =\bbm S^\top &\frac{1}{2} W^\top\\ -S^\top& \frac{1}{2} W^\top \ebm \bbm g\\d \ebm.
 \end{equation}
 System \eqref{eq:blocksystem} is simplified by performing the following block row operations on Row 1 and Row 2 respectively: (1) Multiply Row 1 by $-1$ and subtract Row 2, (2)  multiply Row 2 by -1 and add  Row 1.
 The block matrix system is now
 \begin{equation} \label{eq:blocksyssimp}
     \bbm -\delta_f (x^0;S)-\delta_f(x^0;-S)\\ \delta_f (x^0;S)-\delta_f (x^0;-S)\ebm= \bbm 0&-W^\top \\ 2S^\top &0 \ebm \bbm g\\ d \ebm .
 \end{equation}
 From \eqref{eq:blocksyssimp}, we find 
 \begin{equation*} 
     \frac{1}{2}\left ( \delta_f(x^0;S)-\delta_f (x^0;-S)\right )= S^\top g
 \end{equation*}
 and
 \begin{equation*} 
      \delta_f(x^0;S)+\delta_f (x^0;-S) =W^\top d.
 \end{equation*}
 Now let  $$\delta^c_f(x^0;S)=\frac{1}{2} \bbm f(x^0+s^1)-f(x^0-s^1)\\ \vdots \\ f(x^0+s^k)-f(x^0-s^k)\ebm \in \R^k$$
 and  $$\varepsilon_f(x^0;S)=\bbm f(x^0+s^1)+f(x^0-s^1)-2f(x^0)\\ \vdots \\ f(x^0+s^k)+f(x^0-s^k)-2f(x^0)\ebm \in \R^k.$$
We obtain
 \begin{equation}\label{eq:gradnotisol}
     \frac{1}{2}\left (\delta_f (x^0;S)-\delta_f(x^0;-S)\right )=\delta^c_f(x^0;S)=S^\top g
 \end{equation}
 and
 \begin{equation} \label{eq:hessiannotisol}
      \left (\delta_f(x^0;S)+\delta_f(x^0;-S) \right )=\varepsilon_f(x^0;S)= W^\top d.
 \end{equation}
 
 Analyzing \eqref{eq:gradnotisol}, we see that if $S$ is full row rank, then $(S^\top)^\dagger$ is a left-inverse of $S^\top$  and we can solve \eqref{eq:hessiannotisol} for $g$ by premultiplying both sides by $(S^\top)^\dagger.$ Note that when $S$ is a non-square matrix with full row rank,  then the Moore-Penrose pseudoinverse provides the unique least squares solution of \eqref{eq:gradnotisol} \cite[P. 453]{Horn1990}. Similarly, if $W$ is full row rank then $(W^\top)^\dagger$ is a left-inverse of $W^\top$ and we can solve \eqref{eq:hessiannotisol} for $d$ by premultiplying both sides of \eqref{eq:hessiannotisol} by $(W^\top)^\dagger.$ 
 We are now ready to introduce an explicit  formula to approximate the gradient and a formula to approximate the diagonal of the Hessian based on \eqref{eq:gradnotisol} and \eqref{eq:hessiannotisol}.
\begin{df}[Generalized centered simplex gradient and the centered simplex Hessian diagonal]\label{def:cshd}
Let $f:\dom f \subseteq \R^n\to\R,$ $x^0 \in \dom f$ be the point of interest, $S=\bbm s^1&s^2&\cdots&s^k \ebm$ be in  $\R^{n \times k}$ and $W=\bbm s^1 \odot s^1&\cdots& s^k \odot s^k  \ebm \in \R^{n \times k}.$  Assume that $x^0 \pm s^i  \in \dom f$ for all $i$. 
 The generalized centered simplex gradient of $f$ at $x^0$ over $S$ is denoted by $\nabla_c f(x^0;S)$ and defined by
\begin{equation*}
\nabla_c f(x^0;S)=(S^\top)^\dagger\delta^c_f(x^0;S) \in \R^n.
\end{equation*}
 The centered simplex Hessian diagonal  of $f$ at $x^0$ over $S$,  denoted $\cshd$ is  a vector in $\R^n$ given by 
\begin{align*}
    \cshd&=(W^\top)^\dagger \varepsilon_f(x^0;S).
\end{align*}
\end{df}
Note that the CSHD  can be obtained for only one additional function evaluation when the GCSG has been already computed (by evaluating $f$ at $x^0$). If the function value at the point of interest $x^0$ is known, then the CSHD can be obtained for free in terms of function evaluations. 

The formulae in the previous  definition involves the Moore-Penrose pseudinverse of an $n \times k$ real matrix. A computationally inexpensive approach to compute the Moore-Penrose pseudoinverse is to use  the singular value decomposition of the matrix \cite{wang2018}. For instance, one of the most advanced method, the \emph{Golub-Reisch SVD} algorithm,  requires $4n^2k+8nk^2+9k^3$ flops \cite{Golub1996}. Therefore,  the previous formulae for the  GCSG and CSHD cannot be computed directly in $\mathcal{O}(n)$ flops. The reader is refer to \cite{coope2021gradient}  for more details  on how to rewrite the  formulae so that they can be computed in $\mathcal{O}(n)$ flops when $S$ has a specific form.

Next, we introduce an error bound for the CSHD. Analyzing the error bound, we will get a better understanding  of the matrix $S$ that should be used  to obtain an accurate approximation  of the diagonal entries of the Hessian. In the following theorem, the parameter $\Delta_S$ is the radius of the matrix $S=\bbm s^1&\cdots&s^k \ebm \in \R^{n \times k}.$ That is $$\Delta_S= \max_{i \, \in \{1, 2, \dots k\}}\Vert s^i \Vert.$$
\begin{thm}[Error bound for the CSHD] \label{thm:ebcsdh}
Let $f:\dom f\subseteq \R^n \to \R$ be $\mathcal{C}^4$ on an open domain containing $B(x^0;\Delta_S)$ where $x^0 \in \dom f$ is the point of interest and $\Delta_S>0$ is the radius of $S=\bbm s^1&\cdots s^k\ebm \in \R^{n \times k}.$  Let $W=S \odot S \in \R^{n \times k}.$ Denote by $L_{\nabla^3 f} \geq 0$ the Lipschitz constant of $\nabla^3 f$ on $B(x^0;\Delta_S)$. If $W$ has full row rank, then
\begin{align}\label{eq:EB}
   \left \Vert \cshd-\diagth \right \Vert \leq \left \Vert (\widetilde{W}^\top)^\dagger \right  \Vert \left ( \frac{k}{12}  L_{\nabla^3 f} \Delta_S + 2 \sum_{i=1}^k \left \vert (\hat{s}^i)^\top U[\nabla^2 f(x^0)] \hat{s}^i \right \vert \right )
\end{align}
where $\widetilde{W}=W/\Delta_S^2$ and $\hat{s}^i=s^i/\Delta_S$ for all $i \in \{1, 2, \dots, k\}.$
\end{thm}
\begin{proof}
Since  $W$ has full row rank, $W^\top$ has full column rank and so  $(W^\top)^\dagger$ is a left-inverse of $W^\top.$ We obtain
\begin{align*}
    \left \Vert \cshd-\diagth \right  \Vert&= \left \Vert (W^\top)^\dagger \varepsilon_f(x^0;S)-\diagth \right \Vert\\
    &=\left \Vert (W^\top)^\dagger  \right \Vert \left \Vert \varepsilon_f(x^0;S)- W^\top \diagth \right \Vert.
\end{align*}
Let us investigate the vector $\varepsilon_f(x^0;S)-\diagth \in \R^k$. Row $i$  of this vector can be written as $f(x^0+s^i)+f(x^0-s^i)-2f(x^0)-(w^i)^\top \diagth.$ By Taylor's Theorem, we  know 
\begin{align}\label{eq:fx0psi}
    f(x^0+s^i)=f(x^0)+\nabla f(x^0)^\top s^i+ \frac{1}{2} (s^i)^\top \nabla^2 f(x^0) s^i +(s^i)^\top \nabla^3 f(x^0)[s^i]s^i +R_3(x^0+s^i)
\end{align}
where $R_3(x^0+s^i)$ is the remainder term defined as in \cite[Theorem 1.14]{Burden2016} and $\nabla^3 f(x^0)[s^i] \in \R^{n \times n}$ is a directional Hessian. That is 
\begin{align*}
    \nabla^3 f(x^0)[s^i]&=\lim_{\tau \to 0} \frac{\nabla^2 f(x^0+\tau s^i)-\nabla^2 f(x^0)}{\tau}.
\end{align*}
Similarly, by Taylor's Theorem we may write
\begin{align}\label{eq:fx0msi}
    f(x^0-s^i)=f(x^0)-\nabla f(x^0)^\top s^i+ \frac{1}{2} (s^i)^\top \nabla^2 f(x^0) s^i -(s^i)^\top \nabla^3 f(x^0)[s^i]s^i+ R_3(x^0-s^i).
\end{align}
Adding \eqref{eq:fx0psi} and \eqref{eq:fx0msi} together and rearranging, we obtain
\begin{align*}
    f(x^0+s^i)+f(x^0-s^i)-2f(x^0)&=(s^i)^\top \nabla^2 f(x^0) s^i +R_3(x^0+s^i)+R_3(x^0-s^i).
\end{align*}
Now,  note that $\nabla^2 f(x^0)=\Diagth+\offdiag.$ We get
\begin{align*}
    &f(x^0+s^i)+f(x^0-s^i)-2 f(x^0)-(s^i)^\top \Diagth s^i \\
    &\quad \quad \quad = (s^i)^\top \offdiag s^i + R_3(x^0+s^i)+R_3(x^0-s^i).
\end{align*}
The matrix $\offdiag$ can be written as $\offdiag=\upperh+\upperh^\top.$ Hence, we have $(s^i)^\top \offdiag s^i= 2 (s^i)^\top \upperh s^i$ for all $i$. Therefore, we obtain the equation
\begin{align*}
    &f(x^0+s^i)+f(x^0-s^i)-2 f(x^0)-(w^i)^\top \diagth \\
    &\quad \quad=2(s^i)^\top \upperh s^i + R_3(x^0+s^i)+R_3(x^0-s^i).
\end{align*}
It follows that
\begin{align*}
   &\left  \Vert \varepsilon_f(x^0;S)- W^\top \diagth \right \Vert \\
   &\quad=\sqrt{\sum_{i=1}^k  \left ( 2(s^i)^\top \upperh s^i + R_3(x^0+s^i)+R_3(x^0-s^i)\right )^2 } \\
    & \quad \leq \sum_{i=1}^k \vert 2(s^i)^\top \upperh s^i + R_3(x^0+s^i)+R_3(x^0-s^i) \vert \\
    & \quad\leq  \sum_{i=1}^k \vert 2(s^i)^\top \upperh s^i\vert + \sum_{i=1}^k \vert R_3(x^0+s^i)+R_3(x^0-s^i) \vert \\
    & \quad\leq 2 \sum_{i=1}^k \vert (s^i)^\top \upperh s^i \vert + \sum_{i=1}^k \frac{L_{\nabla^3 f}}{12}\Delta_S^4 \\
    & \quad =2 \sum_{i=1}^k \vert (s^i)^\top \upperh s^i \vert +\frac{k}{12} L_{\nabla^3 f} \Delta_S^4.
\end{align*}
Finally,  the absolute error is 
\begin{align*}
    \left \Vert \cshd -\diagth \right \Vert&\leq \left \Vert (W^\top)^\dagger \right  \Vert \left ( 2 \sum_{i=1}^k \vert (s^i)^\top \upperh s^i \vert + \frac{k}{12} L_{\nabla^3 f} \Delta_S^4 \right )\\
    &= \left \Vert (\widetilde{W}^\top)^\dagger  \right \Vert \left ( 2 \sum_{i=1}^k \vert (\hat{s}^i)^\top \upperh \hat{s}^i \vert + \frac{k}{12} L_{\nabla^3 f} \Delta_S^2 \right ) \qedhere
\end{align*}
\end{proof}
Note that the error bound given in \eqref{eq:EB} does not necessarily need to go to zero as the radius $\Delta_S$ tends to zero. Hence, the true absolute error does not need  to go to zero as $\Delta_S \to 0$. To ensure that the true absolute error goes to zero  when $\Delta_S \to 0,$ we need  the term $\sum_{i=1}^k \vert (\hat{s}^i)^\top \upperh \hat{s}^i \vert$ to vanish. Let us introduce a type of matrix that guarantees that the term $\sum_{i=1}^k \vert (\hat{s}^i)^\top \upperh \hat{s}^i \vert=0.$
\begin{df}[Lonely matrix] 
The matrix $S \in \R^{n \times k}$ is a lonely matrix if there is exactly one nonzero entry in each column of $S.$
\end{df}
\begin{lem} \label{lem:strictupper}
Let $S=\bbm s^1&s^2&\cdots&s^k \ebm$ be a lonely matrix in $\R^{n \times k}.$ Let $U \in \R^{n \times n}$ be a strictly upper triangular  matrix. Then
\begin{align*}
    \sum_{i=1}^k \vert (s^i)^\top U s^i \vert &=0.
\end{align*}
\end{lem}
\begin{lem} \label{lem:Wfullrowrank}
Let $S=\bbm s^1&s^2&\cdots&s^k \ebm$ be a lonely matrix in $\R^{n \times k}$ with full row rank. Then  $W=S \odot S$ has full row rank.
\end{lem}
Obviously, lemma \ref{lem:strictupper}  applies to the set $S=\Id_n$ as it is a diagonal matrix with nonzero entries. Lemma \ref{lem:Wfullrowrank} shows that it is easy to guarantee that $W$ is full row rank by simply taking $S$ to be a lonely matrix with full row rank.   The following corollary provides an error bound when $S$ is a lonely matrix with full row rank.
\begin{cor}\label{cor:lonmatrix}
Let $f:\dom f\subseteq \R^n \to \R$ be $\mathcal{C}^4$ on an open domain containing $B(x^0;\Delta_S)$ where  $x^0 \in \dom f$ is the point of interest and $\Delta_S>0$ is the radius of the matrix $S=\bbm s^1&s^2&\cdots&s^k \ebm \in \R^{n \times k}$. Let $W=S \odot S \in \R^{n \times k}.$ Denote by $L_{\nabla^3 f}\geq 0$ the Lipschitz constant of $\nabla^3 f$ on $B(x^0;\Delta_S)$.  If $S$ is a lonely matrix with full row rank, then
\begin{align*}
    \left \Vert \cshd-\diagth \right \Vert \leq  \left  \Vert (\widetilde{W}^\top)^\dagger \right \Vert \frac{\sqrt{k}}{12}  L_{\nabla^3 f} \Delta_S^2.
\end{align*}
\end{cor}
The previous error bound shows that the CSHD is a $\mathcal{O}(\Delta_S^2)$ accurate approximation of the diagonal entries of the Hessian. Under these assumptions, the optimizer can control the accuracy of the approximation technique. Indeed, as $\Delta_S$ tends to zero, the error bound goes to zero. Hence, the true absolute error needs to go to zero as $\Delta_S \to 0.$ Since the error bound involves the lipschitz constant $L_{\nabla^3 f},$ the CSHD is perfectly accurate when $f$ is a polynomial of  degree less than 4.

In the next section, numerical examples are provided. We will use different matrices $S$  and different points of interest to approximate the CSHD. The relative errors obtained for each $S$ and each point of interest  will be compared.

\section{Numerical examples} \label{sec:num}

In this section,  we begin by conducting a numerical experiment  on the Rosenbrock function $f:\R^2 \to \R:(y_1,y_2) \mapsto (1-y_1)^2+100(y_2-y_1^2)^2$. The CSHD is computed to approximate the diagonal of the Hessian at two points: $x^1=\bbm 1.1&1.1^2+10^{-5}\ebm^\top,$ and $x^2=\bbm 0.9&0.81 \ebm^\top$. Four different matrices $S$ are employed:
\begin{enumerate}
    \item The coordinate basis  in $\R^{n \times n}$ ($\CB$): $\CB=\Id_n.$
    \item A regular basis in $\R^{n \times n}$ ($\RB$): $\RB= \sqrt{\frac{n+1}{n}}\left (\Id_n-\frac{1}{n} \left (1-\sqrt{\frac{1}{n+1}} \right ) \one_n \one_n^\top \right ).$
    \item A coordinate minimal positive basis in $\R^{n \times n+1}$ ($\CMPB$):  $\CMPB=\bbm \Id_n &-\one_n\ebm.$ 
    \item A regular minimal positive basis in $\R^{n \times n+1}$ ($\RMPB$):  $\RMPB=\bbm \RB &-\RB \one_n \ebm.$
\end{enumerate}
Note that the previous four matrices are full row rank, but  only $\CB$ is a lonely  matrix. 
 First, the point of interest $x^1$ is considered and the positive parameter $h^1$ is set to $10^{-3}$. The CSHD is computed using four different sets $S$:  $h^1\CB, h^1\RB, h^1\CMPB,$ and $h^1\RMPB.$ The experiment is repeated with  a new point of interest, $x^2=\bbm 0.9&0.81\ebm^\top.$ The  positive parameter $h^2$ is set to $10^{-6}.$  To compare the accuracy obtained from each set $S$ and each point of interest in $\{x^1, x^2\},$ the relative error will be used instead of the absolute error. This choice is justified by the fact that the true value of the diagonal entries of the Hessian are large in this example. For a given point of interest $x^0$ and matrix $S$, the relative error for a function $f$ at a point $x^0$ over $S$ is defined as 
\begin{align*}
\RER f(x^0;S)&= \frac{\left \Vert \cshd -\diagth \right \Vert}{ \left \Vert \diagth \right  \Vert }.
\end{align*}
Table \ref{table:RE} provides the results obtained.
\begin{table}[ht] 
\caption{\textbf{Relative Error of the diagonal entries of the Hessian}}
\center{
\begin{tabular}{|p{2.3cm}|p{3.0cm}|p{3.0cm}|}
 \hline
\textbf{Set $T$} &\textbf{$\RER f(x^1;h^1T)$}&\textbf{$\RER f (x^2;h^2T)$}\\
 \hline
 CB&$2.02 \times 10^{-7}$&$1.18 \times 10^{-9}$\\
 \hline
 RB&$3.14 \times 10^{-1}$&$3.74 \times 10^{-1}$\\
 \hline
 CMPB&$4.19 \times 10^{-1}$&$4.99 \times 10^{-1}$\\
 \hline
 RMPB&$1.78 \times 10^{-7}$&$3.39 \times 10^{-9}$\\
 \hline
\end{tabular}}\label{table:RE}
\end{table}
Analyzing Table \ref{table:RE}, we observe that the sets  $\RB$ and $\CMPB$ provide poor approximations of the diagonal entries of the Hessian. Let us check if it is possible to decrease the relative error  by decreasing the positive parameter $h$. Table \ref{table:RElim} provides the limit of the relative error as $h \to 0$ and  the infimum of the relative error.
 Maple 2020 is employed to conduct this experiment. In Table \ref{table:RElim}, $\POI$ stands for point of interest. Table \ref{table:RElim}   shows that the limit as $h$ tends to zero is not  equal to 0 when RB, CMPB and RMPB are utilized.  It is also interesting to note that the minimum of $\RER f(x^1;h \RMPB)$ is approximately $5.71 \times 10^{-10}$ but the limit when $h$ tends to zero of $\RER f(x^1;h\RMPB)$ is infinity (this is not due to numerical errors).
\begin{table}[ht] 
\caption{\textbf{Comparing four matrices $S$ to approximate  the diagonal of the Hessian}}
\center{
\begin{tabular}{|p{3.0cm}|p{3.6cm}|p{3.1cm}|}
 \hline
\textbf{$(\POI, S)$} & $\lim_{h \to 0} \RER f(\POI;S)$&$\inf_{h} \RER f(\POI;S)$\\
 \hline
 $(x^1;h\CB)$&0 &0  \\
 \hline
  $(x^2;h\CB)$&0 &0  \\
 \hline
  $(x^1;h\RB)$&$3.14 \times 10^{-1}$ &$3.14 \times 10^{-1}$  \\
 \hline
  $(x^2;h\RB)$&$3.74 \times 10^{-1}$ &$3.74 \times 10^{-1}$  \\
 \hline
  $(x^1;h\CMPB)$&$4.19 \times 10^{-1}$ &$2.96 \times 10^{-1}$  \\
 \hline
  $(x^2;h\CMPB)$&$5.00 \times 10^{-1}$ &$3.53 \times 10^{-1}$ \\
 \hline
  $(x^1;h\RMPB)$&$+\infty $&$5.71 \times 10^{-10}$ \\
 \hline
  $(x^2;h\RMPB)$&$4.65 \times 10^{-10}$ & $4.65 \times 10^{-10}$\\
 \hline
\end{tabular}}\label{table:REinv} \label{table:RElim}
\end{table}

The results of Table \ref{table:RElim}  agree with the error bound defined \eqref{eq:EB}.  For instance, consider the error bound for $\diag [\nabla_c^2 f(x^0;h \CMPB)].$ First, note that the radius of the set $S=h\cdot \CMPB$ is $\Delta_S=h\sqrt{2}.$ Looking at the term $\sum_{i=1}^{n+1} \vert (\hat{s}^i)^\top \upperh \hat{s}^i \vert$ in \eqref{eq:EB}, we see that 
 \begin{align*}
     \sum_{i=1}^{3} \vert (\hat{s}^i)^\top \upperh \hat{s}^i \vert&=h^2\vert \widehat{\one_2}^\top \upperh \widehat{\one_2} \vert\\
     &=\frac{1}{2} \left \vert \one_2^\top \upperh \one_2  \right \vert\\
     &=\frac{1}{2}  \left \vert \bbm 1&1\ebm^\top \bbm0&-440\\0&0\ebm \bbm 1\\1\ebm \right \vert=220.
 \end{align*}
 Observe that, independently of the value $h$, the sum $\sum_{i=1}^{3} \vert (\hat{s}^i)^\top \upperh \hat{s}^i \vert$ is always equal to 220 when the matrix CMPB is used on this specific function $f$ at the point of interest $x^1$.
 Therefore, the limit of the error bound in \eqref{eq:EB} as $h$ tends to zero (and so $\Delta_S$ tends to zero) is greater than zero.  Hence, the error bound does not go to  zero as $h$ tends to zero and so nothing guarantees that the true absolute error goes to zero as $h \to 0$.

The previous numerical results and the error bound developed in \eqref{eq:EB} strongly suggests  that $\CB$, or any lonely matrix with full row rank,  should be used  when  approximating the diagonal entries of the Hessian via the technique $\cshd$. Based on the results obtained in Table \ref{table:RE} and Table \ref{table:RElim}, the matrix RMPB seems to provide better results than $\CMPB$ and $\RB$. The following example will show that RMPB can behave badly.
\begin{ex} \label{ex:RMPB}
Let $f:\R^3 \to \R:(y_1,y_2,y_3) \mapsto e^{y_1y_2y_3}.$ Let $x^0=\bbm 3&2&1 \ebm^\top$ be the point of interest. Let $h>0$ be the shrinking parameter and let the matrix $S \in \R^{3 \times 4}$ be 
\begin{align*}
S&=h\RMPB=h\bbm \frac{5\sqrt{3}}{9} &-\frac{\sqrt{3}}{9}&-\frac{\sqrt{3}}{9}&-\frac{\sqrt{3}}{3}\\-\frac{\sqrt{3}}{9}&\frac{5\sqrt{3}}{9}&-\frac{\sqrt{3}}{9}&-\frac{\sqrt{3}}{3}\\-\frac{\sqrt{3}}{9}&-\frac{\sqrt{3}}{9}&\frac{5\sqrt{3}}{9}&-\frac{\sqrt{3}}{3}\ebm.
\end{align*}
Then
\begin{align*}
    \lim_{h \to 0} \RER f(x^0;h \RMPB)&\approx 1.33 \times 10^{-1},
\end{align*}
and
\begin{align*}
    \min_h \RER f(x^0;h\RMPB)&\approx 1.30 \times 10^{-1}. 
\end{align*}
The minimal relative error $\RER f(x^0;h\RMPB)$  is attained at $h^*\approx 0.0883.$ Table \ref{table:diffh} compares the relative error $\RER f(x^0;h \RMPB)$ and $\RER f(x^0;h \CB)$ for different values of $h$.
\begin{table}[ht]
\caption{\textbf{Comparing RMPB and CB as $h$ tends to zero}}
\center{
\begin{tabular}{|p{2.0cm}|p{3.5cm}|p{3.0cm}|}
 \hline
\textbf{h} & $\RER f(x^0;h\RMPB)$& $\RER f(x^0;h \CB)$\\
 \hline
 1&$5.93\times 10^{1}$& $9.79 \times 10^0$\\
 \hline
 $10^{-1}$& $1.31 \times 10^{-1}$&$2.93 \times 10^{-2}$\\
 \hline
 $10^{-2}$&$1.33 \times 10^{-1}$ & $2.90 \times 10^{-4}$\\ 
 \hline
 $10^{-3}$& $1.33 \times 10^{-1}$& $2.90 \times 10^{-6}$\\ 
 \hline
 $10^{-4}$&$1.33 \times 10^{-1}$ & $2.95 \times 10^{-8}$\\ 
 \hline
 $10^{-\infty}$&$1.33 \times 10^{-1}$&0\\
 \hline
\end{tabular}}\label{table:diffh}
\end{table}
\end{ex}
This example agrees with the error bound introduced in \eqref{eq:EB} as nothing guarantees that the term $\sum_{i=1}^{n+1} \vert (\hat{s}^i)^\top \upperh \hat{s}^i \vert$ is equal to zero when the matrix $\RMPB$ is used to form the sample set of points.
\section{Conclusion} \label{sec:conc}
Based on the error bound developed in Theorem \ref{thm:ebcsdh} and our numerical results, the coordinate basis is the best choice of the  four sets tested to approximate the diagonal entries of the Hessian via the technique $\cshd.$  More generally, a set of directions that forms a lonely matrix with full row rank should be used to approximate the diagonal entries of the Hessian. When this is the case, the CSHD is $\mathcal{O}(\Delta_S^2)$ accurate.

It is worth mentioning that approximating the diagonal of the Hessian and not the entire Hessian may be misleading  when  the Hessian of $f$ is not diagonally dominate. For example, consider the function $f(y_1,y_2)=\alpha y_1 y_2$ where $\alpha$ is a nonzero scalar. Then $$\nabla^2 f(y_1,y_2)=\bbm 0 &\alpha\\ \alpha &0 \ebm.$$ In this case, the diagonal of the Hessian does not provide any information about the function $f$ even though $f$ is not equal to the zero function. 

Future research directions could include testing a  diagonal quadratic model  in a  model-based trust region method \cite[Chapter 11]{audet2017derivative} and  comparing its performance to  a (full) quadratic model. It may be valuable to compare  these two approaches in terms of function evaluations necessary to obtain a solution within a certain accuracy.
\section*{Acknowledgement}

Jarry-Bolduc's research is partially funded by the Natural Sciences and Engineering Research Council (NSERC) of Canada, Discover Grant \#2018-03865. Jarry-Bolduc would like to acknowledge UBC for the funding received through the University Graduate Fellowship award.
\bibliographystyle{plain}
\bibliography{Bibliography}

\def\cprime{$'$} \def\cprime{$'$}
\begin{thebibliography}{10}

\bibitem{Audet2014}
C.~Audet.
\newblock A survey on direct search methods for blackbox optimization and their
  applications.
\newblock In {\em Mathematics Without Boundaries}, pages 31--56. Springer,
  2014.

\bibitem{Audet2017}
C.~Audet and W.~Hare.
\newblock {\em Derivative-free and blackbox optimization}.
\newblock Springer Series in Operations Research and Financial Engineering.
  Springer, Switzerland, 2017.

\bibitem{audet2017derivative}
C.~Audet and W.~Hare.
\newblock {\em Derivative-free and {B}lackbox {O}ptimization}.
\newblock Springer, 2017.

\bibitem{Billups2013}
S.~Billups, J.~Larson, and P.~Graf.
\newblock Derivative-free optimization of expensive functions with
  computational error using weighted regression.
\newblock {\em SIAM J. Optim.}, 23(1):27--53, 2013.

\bibitem{Bortz1998}
D.~Bortz and C.~Kelley.
\newblock The simplex gradient and noisy optimization problems.
\newblock In {\em Computational Methods for Optimal Design and Control}, pages
  77--90. Springer, 1998.

\bibitem{Burden2016}
A.~Burden, R.~Burden, and J.~Faires.
\newblock Numerical analysis 10/e ie.
\newblock {\em Brooks/Cole Cengage Learning}, 2016.

\bibitem{Conn22008}
A.~Conn, K.~Scheinberg, and L.~Vicente.
\newblock Geometry of sample sets in derivative-free optimization: polynomial
  regression and underdetermined interpolation.
\newblock {\em IMA J. Num. Anal.}, 28(4):721--748, 2008.

\bibitem{Conn2009}
A.~Conn, K.~Scheinberg, and L.~Vicente.
\newblock {\em Introduction to Derivative-free Optimization}, volume~8.
\newblock Siam, 2009.

\bibitem{coope2021gradient}
I.~Coope and R.~Tappenden.
\newblock Gradient and diagonal hessian approximations using quadratic
  interpolation models and aligned regular bases.
\newblock {\em Numerical Algorithms}, pages 1--25, 2021.

\bibitem{Custodio2017}
A.~Cust{\'o}dio, K.~Scheinberg, and L.~Vicente.
\newblock Methodologies and software for derivative-free optimization.
\newblock {\em Advances and Trends in Optimization with Engineering
  Applications}, pages 495--506, 2017.

\bibitem{Golub1996}
G.~Golub and C.~Van~Loan.
\newblock {\em Matrix Computations (3rd Ed.)}.
\newblock Johns Hopkins University Press, USA, 1996.

\bibitem{hare2020}
W.~Hare and G.~Jarry-Bolduc.
\newblock Calculus identities for generalized simplex gradients: Rules and
  applications.
\newblock {\em SIAM J. Optim.}, 30(1):853--884, 2020.

\bibitem{hare2020error}
W.~Hare, G.~Jarry-Bolduc, and C.~Planiden.
\newblock Error bounds for overdetermined and underdetermined generalized
  centred simplex gradients.
\newblock {\em arXiv preprint arXiv:2006.00742}, 2020.

\bibitem{HareNutiniTesfamariam2013}
W.~Hare, J.~Nutini, and S.~Tesfamariam.
\newblock A survey of non-gradient optimization methods in structural
  engineering.
\newblock {\em Adv. Eng. Soft.}, 59:19--28, 2013.

\bibitem{Horn1990}
R.~Horn and C.~Johnson.
\newblock {\em Matrix analysis}.
\newblock Cambridge university press, 1990.

\bibitem{Kelley1999}
C.~Kelley.
\newblock {\em Iterative Methods for Optimization}, volume~18.
\newblock SIAM, 1999.

\bibitem{larson_menickelly_wild_2019}
J.~Larson, M.~Menickelly, and S.~Wild.
\newblock Derivative-free optimization methods.
\newblock {\em Acta Numer.}, 28:287--404, 2019.

\bibitem{Oeuvray2007}
R.~Oeuvray and M.~Bierlaire.
\newblock Boosters: A derivative-free algorithm based on radial basis
  functions.
\newblock {\em Int. J. Model. Sim.}, 29(1):26--36, 2009.

\bibitem{Powell2004}
M.~Powell.
\newblock Least {F}robenius norm updating of quadratic models that satisfy
  interpolation conditions.
\newblock {\em Math. Program.}, 100(1):183--215, 2004.

\bibitem{Regis2015}
R.~Regis.
\newblock The calculus of simplex gradients.
\newblock {\em Optim. Lett.}, 9(5):845--865, 2015.

\bibitem{Regis2005}
R.~Regis and C.~Shoemaker.
\newblock Constrained global optimization of expensive black box functions
  using radial basis functions.
\newblock {\em J. Global Optim.}, 31:153--171, 2005.

\bibitem{rockwets}
R.~Rockafellar and R.~Wets.
\newblock {\em Variational analysis}.
\newblock Grundlehren der Mathematischen Wissenschaften [Fundamental Principles
  of Mathematical Sciences]. Springer-Verlag, Berlin, 1998.

\bibitem{Schonlau1998}
M.~Schonlau, W.~Welch, and D.~Jones.
\newblock Global versus local search in constrained optimization of computer
  models.
\newblock {\em Lecture Notes Monogr. Ser.}, pages 11--25, 1998.

\bibitem{wang2018}
G.~Wang, Y.~Wei, S.~Qiao, P.~Lin, and Y.~Chen.
\newblock {\em Generalized inverses: theory and computations}, volume~53.
\newblock Springer, 2018.

\bibitem{Wild2011}
S.~Wild and C.~Shoemaker.
\newblock Global convergence of radial basis function trust region
  derivative-free algorithms.
\newblock {\em SIAM J.~Optim.}, 21(3):761--781, 2011.

\end{thebibliography}
\pagebreak

\end{document}